\documentclass[12pt]{article}
\usepackage{amsmath,amssymb,bm,amsthm}
\usepackage{graphicx,xcolor}
\usepackage{array,cite,url}
\usepackage{amsfonts,epstopdf}
\usepackage{booktabs}
\usepackage{array}
\usepackage{caption}
\usepackage[position=top]{subfig}
\usepackage[shortlabels]{enumitem}
\usepackage[bookmarks=false]{hyperref} 

\usepackage[title]{appendix}
 
\usepackage[]{cleveref}
\usepackage{indentfirst}
\usepackage{tikz}
\graphicspath{{./}{./figures/}{./figures/FrameworkSketch/}}

\newcommand{\Ex}{\mathbb{E}}

\newcommand{\rar}{\rightarrow}
\newcommand{\bs}{\boldsymbol}

\newcommand{\U}{\mathcal{U}}

\newcommand{\x}{\mathbf{x}}
\newcommand{\ub}{\mathbf{u}}
\newcommand{\h}{\mathcal{H}}
\newcommand{\R}{\mathbb{R}}

\usepackage{algorithm}
\usepackage[noend]{algpseudocode}
\makeatletter
\def\BState{\State\hskip-\ALG@thistlm}
\makeatother
\usepackage{algorithmicx}

\title{\LARGE \bf
	Generative stochastic modeling of strongly nonlinear flows with non-Gaussian statistics}

\author{Hassan Arbabi and Themistoklis Sapsis
	\thanks{The authors are with the department of Mechanical Engineering,
		Massachusetts Institute of Technology, Cambridge, MA 02139, USA. Correspondence: {\tt\small\{arbabi, sapsis\}@mit.edu}}%
}

\hypersetup{
  colorlinks   = true, 
  urlcolor     = blue, 
  linkcolor    = blue, 
  citecolor   = magenta 
}

\newif\ifcomments

\commentsfalse

\newif\ifchanged
\newcommand\changed[1]{\ifchanged\textcolor{black}{ { #1}}\else\relax\fi}
\changedtrue

\begin{document} 

\maketitle

\begin{abstract}
Strongly nonlinear flows, which commonly arise in geophysical and engineering turbulence, are characterized by persistent and intermittent energy transfer between various spatial and temporal scales. 
These systems are difficult to model and analyze due to combination of high dimensionality and uncertainty, and there has been much interest in obtaining reduced models, in the form of stochastic closures, that can replicate their non-Gaussian statistics in many dimensions. 
Here, we propose a data-driven framework to model stationary chaotic dynamical systems through nonlinear transformations and a set of decoupled stochastic differential equations (SDEs).  Specifically, we use optimal transport to find a transformation from the distribution of time-series data to a multiplicative reference probability measure such as the standard normal distribution. Then we find the set of decoupled SDEs that admit the reference measure as the invariant measure, and also closely match the spectrum of the transformed data. As such, this framework represents the chaotic time series as the evolution of a stochastic system observed through the lens of a nonlinear map. We demonstrate the application of this framework in Lorenz-96 system, a 10-dimensional model of high-Reynolds cavity flow, and reanalysis climate data. 
These examples show that SDE models generated by this framework can reproduce the non-Gaussian statistics of systems with moderate dimensions (e.g. 10 and more), and predict super-Gaussian tails that are not readily available from little training data. These findings suggest that this class of models provide an efficient hypothesis space for learning strongly nonlinear flows from small amounts of data.
\end{abstract}


\clearpage


\section{Data-driven modeling  of dynamical systems}
Scientists and engineers are facing increasingly difficult challenges such as predicting and controlling the earth climate, understanding and emulating the human brain, and design and maintenance of social networks. All these challenges require modeling and analysis of systems that have a large number of degrees of freedom, possibly combined with a considerable amount of uncertainty in parameters.  In doing so, the classical model-based state-space approach for dynamical systems analysis and control falls short due to the computational size of these problems and lack of an accurate model. As a result, the recent decades have seen a great development in the area of  model reduction and data-driven modeling with the focus of making these problems more amenable to computation and analysis.

One of the main challenges in modeling complex systems is centered around strongly nonlinear systems, i.e. chaotic systems with persistent and intermittent energy transfer between multiple scales. Climate dynamics and engineering turbulence give rise to many such problems and have been a major target of model reduction techniques. In particular, using stochastic reduced models has become increasingly popular for these problems because it allows representing turbulent fluctuations,
which arise from nonlinear interactions in many degrees of freedom, in the form of simple systems with random forcing. Stochastic models of climate are proposed  to estimate climate predictability and response to external output \cite{hasselmann1976stochastic,leith1975climate},   explain the observed energy spectrum and transport rates in atmosphere \cite{farrell1994theory,delsole2004stochastic}, facilitate data assimilation \cite{harlim2008filtering} and predict weather patterns \cite{chen2014predicting}. Similar models for analysis of turbulent flows have been suggested in \cite{sreenivasan2002mean,de2007slow,rigas2015diffusive,zare2017colour,khodkar2018data} . The formal analysis in  \cite{majda1999models, majda2001mathematical, majda2006stochastic} has established the validity of stochastic closure models under certain assumptions while emphasizing the need to go beyond the commonly used linear models with additive noise. All the above efforts have shown the great utility and promise of stochastic modeling, but they heavily rely on expert intuition and the underlying model of the system (e.g. linearized Navier-Stokes equations) and therefore their application remains limited to well-understood systems and well-chosen variables.

Data-driven methods provide a shortcut for analysis and control of high-dimensional systems by allowing us to discover and exploit low-dimensional structures that may be masked by the governing equations, or find alternative models with more computational tractability.
A key enabler of data-driven analysis for dynamical systems in the recent years has been the operator-theoretic formalism,\changed{and in particular, the Koopman operator theory  \cite{koopman1931,mezic2005}, which describes the evolution of observables, for computation of  geometric objects and linear predictors from trajectory data. This approach is accompanied by an arsenal of numerical algorithms that extract dynamical information from large data produced by simulations and experiments \cite{schmid2010,rowley2009,williams2015data,noe2013variational,perez2013identification}. 
Due to interpretability and connection with the classical theory of dynamical systems, this methodology has become popular in many fields ranging from fluid dynamics \cite{schmid2011applications} and power networks \cite{raak2015data} to biological pattern extraction \cite{brunton2016extracting} and visual object recognition \cite{erichson2016compressed}. 
The Koopman operator theory is closely related to the Perron-Frobenius operator and Fokker-Planck  equations \cite{LasotaandMackey:1994,dellnitz1999approximation,DasBuch2005} however those viewpoints have enjoyed less widespread use in data-driven applications.  
A central feature of the Koopman operator framework is to use data to find the canonical coordinates which decouple the system dynamics and allow independent and linear representation. This feature is also emphasized in manifold learning methods for data mining which are powered by (diffusion) operator realization and analysis \cite{coifman2006diffusion,singer2008non}. 
However, much of the current operator-theoretic framework relies on the discrete spectrum expansion, and the modeling framework for stationary chaotic systems, where the associated operators possess continuous spectrum, is still in a nascent phase \cite{korda2018data,giannakis2021delay,colbrook2021rigorous}.}

Statistical and machine learning methods offer another large and open-ended avenue for data-driven modeling. Variations of neural networks; including auto-encoders, long short-term memory networks and reservoir computers; are proposed for various tasks like discovering representations of complex systems, closure modeling and short-time prediction of chaotic systems \cite[e.g.]{rico1992discrete,rico1995nonlinear,pathak2018model,wan2018data,lusch2018deep,chattopadhyay2020data}. A big theme of many efforts in this area is to explicitly encode available physical and mathematical information in the structure of learning to reduce the amount of training data and adhere to known constraints \cite{wang2017physics,wan2018machine,raissi2018numerical}.
Furthermore the approximation of operator-theoretic objects from data is naturally connected to the regression problem and there is a growing body of works that use machine learning techniques for approximation of the Koopman or Perron-Frobenius operator \cite{kevrekidis2015kernel,zhao2016analog,takeishi2017learning,yeung2017learning,li2017extended,lusch2018deep,otto2019linearly}. 
All the above studies show the promise of machine learning in modeling dynamical systems but the methods in this category often accompany important caveats such as lack of interpretability and non-robust optimization processes.

In this work, we present a data-driven framework for \emph{generative modeling} of strongly nonlinear flows. In contrast to the most of data-driven modeling frameworks,  which focus on the short-time prediction, our goal is to discover models from data that emulate the statistics of a strongly nonlinear flow in the appropriate sense.
 This problem, which involves matching joint distributions of stochastic processes at arbitrary lags,  is  computationally prohibitive for strongly nonlinear flows. However, inspired by the spectral theory of dynamical systems \cite{mezic2005,korda2018data}, we reformulate this problem as finding a model that replicates the invariant measure as well as the power spectral density of the observed time-series data. The cornerstone  of our approach is using optimal transport of probabilities \cite{parno2018transport,villani2008optimal} to map the data distribution to a multiplicative reference measure and then model each dimension separately using a simple stochastic differential equation (SDE) that mimics the spectrum of the time series. Starting from time-series measurements of observables on a dynamical system, our framework yields a set of randomly forced linear oscillators combined with a nonlinear observation map that produces the same statistics and dynamics as the given time-series data.

Our framework enables a stronger data-driven approach compared to the aforementioned works on stochastic modeling: 
we start with a general space of stochastic models, i.e., linear stochastic oscillators pulled back under polynomial maps, and then select the model that produces closest spectra to the observed data. This offers advantage over the previous works on data-driven modeling of strongly nonlinear systems  \cite{kravtsov2005multilevel,majda2012physics} by treating nonlinearities as a feature of the observation map and not the underlying vector field. As a result, we can use dynamic models with guaranteed stable invariant measures as opposed to quadratic models in \cite{kravtsov2005multilevel}, and avoid using latent variables as opposed to \cite{majda2012physics}. Despite the generality of our framework, our examples show that this choice of hypothesis space is narrow enough to provide a data-efficient and computationally tractable pipeline for learning statistically accurate models of strongly nonlinear flows with tens of dimensions. This is partially due to formulating a major part of the challenge as a probability transport problem in large dimensions which has enjoyed great computational progress in recent years.

We use our framework to compute stochastic models for the Lorenz-96 system, high-Reynolds flow in a 2D cavity and 6-hourly time series of velocity and temperature reanalysis data from the earth atmosphere. \changed{In case of cavity flow, we use the Spectral Proper Orthogonal Decomposition (SPOD) \cite{towne2018spectral} to obtain the 10-dimensional model of the flow. A theoretical contribution of our work is showing a novel connection between this decomposition and spectral decomposition of the Koopman operator for systems with continuous spectrum.}
 Application of our framework to the above examples generates models that closely match the non-Gaussian features of data such as skewness and heavy tails in large dimensions from relatively little data. In case of the cavity flow, the SDE models not only produce the statistics of modal coordinates which are used to train the model, but they also lead to accurate approximation of high-order statistics for pointwise flow measurements. 
The transport-based framework also allows modeling systems that possess heavy tails in their invariant measures. By utilizing this feature in case of the climate data, we are able to extrapolate the heavy tails of the chaotic reanalysis time series.This contribution is specially important since it allows probabilistic characterization of extreme events from short-time observations, which is an outstanding challenge in modeling strongly nonlinear flows \cite{lucarini2016extremes}.


\section{Stochastic generative modeling of strongly nonlinear flows} \label{sec_framework}

\subsection{Problem setup}
Consider a dynamical system given as
\begin{align}\label{eq_dynamics}
\dot{x}&=f(x),\\
y&=g(x) \notag
\end{align}
with the state variable $x$ in the state space $S$ and the output $y$ in $\mathbb{R}^N$.  We use $F^t:S\rar S$ to denote the flow map over time $t$, that is, the mapping that takes the state at time $t_0$ to the state at time $t_0+t$. 
Assume that the trajectories in an open subset of $S$ converge to an attractor $\Omega$ \changed{which supports an invariant and ergodic probability measure denoted by $\mu$, that is, $\mu \big(F^{-t}(B)\big) = \mu(B)$ for all measurable $B\subset \Omega$, $\mu(\Omega)=1$, and $\mu(\Omega-C)=0$ or 1 if $C$ in an invariant subset of $\Omega$. }

We assume that instead of $x$, we can only measure  $y$ which is a random variable with distribution $\nu=g_\#\mu$, that is, the pushforward of $\mu$ under the \changed{(measurable)} observation map $g$, defined as
\begin{align}
\nu(g<a)=\mu\bigg(g^{-1}\big((-\infty,a)\big)\bigg)
\end{align}
We let $U^t$ be the Koopman operator of this system which describes the evolution of random variables with the dynamics, that is, $U^tg:=g\circ F^t$ \cite{koopman1931,mezic2005}. \changed{Moreover, we assume the observables we study are  square integrable with respect to the invariant measure $\mu$ so that we can utilize the spectral results developed for the Koopman operator on such Hilbert spaces.}    The collection of random variables $\{U^t g\}_{t\in\mathbb{R}}$ is a stationary stochastic process and time-series of $y$ is a realization of this process \cite[e.g.]{doob1953stochastic}.

Given the time series of $y$, our goal is to construct a dynamical system that generates a statistically similar stochastic process. 
A \emph{perfect generative model} for this system will generate a stochastic process that has the same joint distributions as  this process for arbitrary combination of lags.  \changed{To be more precise, let 
\begin{equation}\label{eq_process}
P_{\boldsymbol{\tau},\mathbf{a}}[Z^t]= P(Z^0>\mathbf{a}_0,Z^{\tau_1}>\mathbf{a}_1,\ldots,Z^{\tau_n}>\mathbf{a}_n)
\end{equation}
be the joint distribution of the stochastic process $\{Z^t\}_{t\in\R}$, with $\boldsymbol{\tau} =(\tau_1,\ldots,\tau_n)\in \R^{n}$,  $\mathbf{a}_i \in \R^{N}$ for $i=0,\ldots,n$ and $Z^{\tau_i}>\mathbf{a}_i$ denotes element-wise inequality between the vectors $\mathbf{a}_i$ and $Z^{\tau_i}$}.
Then $\{Z^t\}_{t\in\R}$ is a perfect model for our time series  if
\begin{equation}\label{eq_process}
P_{\boldsymbol{\tau},\mathbf{a}}[U^tg]=  P_{\boldsymbol{\tau},\mathbf{a}}[Z^t]
\end{equation}
for all choices of $n$, $\boldsymbol{\tau}$ and $\mathbf{a}$. Searching for a good model of time series using this characterization is computationally formidable since it requires presentation and matching of many, and in principle infinite number of, joint distributions.

In this work, we appeal to another characterization of stochastic processes which is motivated by spectral analysis of dynamical systems. Recall that the power spectral density (PSD) of a stationary process is given by Fourier transform of its covariance, e.g., for the above process we have\footnote{In this paper, for simplicity of presentation, we have assumed that the power spectral measure is absolutely continuous and hence representable by a density. Same theoretical arguments hold if we use the power spectral measure. }
\begin{equation}\label{eq_PSD}
S_g(\omega) = \int_{-\infty}^{\infty}e^{-i\omega \tau} \text{Cov}(g,U^\tau g) d\tau.
\end{equation}
It turns out that given the knowledge of the the underlying attractor $\Omega$ and the invariant measure $\mu$, the PSD of observable $g$ contains all the information about the dynamics of the system that can be gleaned from observing $g$. More precisely, assume $g$ is square integrable with respect to $\mu$ and then define 
\begin{equation}\label{eq_Hf}
\h_g=\overline{\text{span} \{g,U^{\tau}g,U^{-\tau}g,U^{2\tau}g,U^{-2\tau}g,\ldots\}}.
\end{equation}
where $\tau$ is the sampling interval.
In other words, $\h_g$ is the closure of random variables that can be represented as a linear combination of $g$ and its history. Then with the knowledge of the PSD of $g$, we can uniquely determine the Koopman operator $U$ restricted to $\h_g$ {\cite[Proposition 1]{korda2018data}}, implying that we can predict the time evolution of every random variable in $\h_g$. In case that $g$ is informative enough \changed{(or more precisely, $*$-cyclic under the action of $U^t$ \cite{maccluer2008elementary})},  $\h_g$ contains all square-integrable random variables including the state variables.

Guided by the above observation, we reformulate the generative modeling for a stationary stochastic processes as follows: we seek a model that replicates the same 1st-order distribution (i.e. $\nu$, the invariant measure observed through the lens of $g$), and the same PSD as the time series of $g$. We are interested in finding such models for strongly nonlinear flows with many dimensions. 
For computational tractability in these problems, we propose a hypothesis space which consists of linear stochastic oscillators observed under the inverse of (invertible) multivariate polynomial maps. To find the best model within this space, we first use optimal transport to discover the appropriate observation map that generates the distribution of the target stochastic process. Then, we optimize the parameters of the oscillator to best match the PSD of the target stochastic process. 
Our results and analysis show that this choice of hypothesis space and learning process is data efficient and computationally economic for learning strongly nonlinear flows.



\subsection{Modeling with optimal transport and spectral matching}\label{sec_TMmodeling}
The first step in our framework is to find the mapping $T:\mathbb{R}^N \rar \mathbb{R}^N$ that takes the random variable $y\sim \nu$ to another random variable, denoted as $q$, which has a standard normal distribution, that is, 
\begin{equation}\label{eq_Tq}
q=T(y)\sim \pi,
\end{equation}
where $\pi$ is the standard normal distribution.
 In other words, we seek the change of variables that makes the attractor look like a normalized Gaussian distribution in each direction.  We use the theory of \emph{optimal transport} to find this change of variable. 
This theory studies the choice of maps which minimize some notion of cost for carrying a probability measure to another \cite{villani2008optimal,peyre2019computational}. The existence of such a map generally depends on the regularity of the measure $\nu$. We assume that $\nu$ is absolutely continuous with respect to the Lebesgue measure on $\mathbb{R}^N$  and possesses finite second-order moments which guarantees the existence of the transport map  \cite{villani2008optimal}.  
For numerical approximation of this map, we use the computational approach developed in \cite{parno2018transport,marzouk2016sampling} (also see \cite{el2012bayesian,spantini2018inference,morrison2017beyond}). This framework leverages the structure of a specific class of optimal transport maps, known as Knothe--Rosenblatt rearrangement, which leads to advantageous computational properties.  The map is assumed to be a lower triangular, differentiable and monotone function given as
 \begin{align}\label{eq_triangular}
T(y_{1},y_{2},\ldots,y_{N})=
\begin{bmatrix}
T_1(y_1) 
\\ T_2(y_1,y_2)
\\ \vdots
\\T_N(y_{1},y_2,\ldots,y_{N})
\end{bmatrix},
\end{align}
where each component is approximated with a multivariate polynomial of a prescribed degree. To ensure the monotonicity, and hence invertibility of $T$, we use the integrated squared parameterization, i.e.,
 \begin{align}\label{eq_Tcomp}
T_i(y_{1},y_2,\ldots,y_{i}) = c(\mathbf{a}_c;y_1,\ldots,y_{i-1})+\int_{0}^{y_i}\big(h(\mathbf{a}_h;y_1,\ldots,y_{i-1},t)\big)^2dt
\end{align}
with 
 \begin{align}\label{eq_Tcomp}
c =  \boldsymbol{\Phi}_c(y_1,\ldots,y_{i-1}) \mathbf{a}_c,\quad h =  \boldsymbol{\Phi}_e(y_1,\ldots,y_{i-1}) \mathbf{a}_h,
\end{align}
where basis functions in $\boldsymbol{\Phi}_c$ are multivariate Hermite polynomials, and $\boldsymbol{\Phi}_e$ are Hermite functions extended with constant basis functions \cite{peherstorfer2019transport,transportmaps_web}.

The parameters of $T$, i.e. expansion coefficients $\mathbf{a}_c$ and $\mathbf{a}_e$,  are determined by maximum likelihood estimation: let $\tilde T$ be a candidate for $T$, and consider the pullback of $\pi$ under this map,
\begin{align}\label{eq_pullback}
\tilde \nu(y)=\pi(\tilde T(y))|\det\nabla \tilde T(y)|.
\end{align}
Then the map $T$ is found by minimizing the Kullback-Leibler divergence between this pullback and the data distribution $\nu$, i.e.,
 \begin{align}\label{eq_TMobj}
T=\arg\min_{\tilde T}    \Ex_{\nu}\bigg[\log \frac{\nu(y)}{\tilde \nu(y)} \bigg]=\arg\min_{\tilde T}    \Ex_{\nu}\bigg[-\log{\tilde \nu(y)} \bigg]
\end{align}
where we dropped $\nu(y)$ in the second equality since it is independent of optimizer $\tilde T$.
The expectation with respect to $\nu$ is approximated using the time average of data due to ergodicity.
The key feature of this computational setup is that it leads to an optimization problem which can be solved for each $T_i$ separately \cite{marzouk2016sampling,parno2018transport}, thereby allowing an efficient computation for large dimensions and long time series.

After we have found $T$, we fit a system of stochastic oscillators to the time series of the random variable \changed{ $q=T(y)$}. Since we chose $\pi$ to be a multiplicative measure, we can do this fitting for each component of $q$ independently. \changed{The independent modeling of each coordinate is a key aspect of this work which enables computational treatment of high-dimensional systems, as opposed to direct discovery of a dynamical systems which would generate the same high-dimensional joint distribution as the data.} We consider the systems of forced linear oscillators,
\begin{align}\label{eq_SDE}
\ddot{q}_{j}+ \beta_j \dot{q}_j +k_j q_j= \sqrt{2D_j} \dot{W}_j, \quad \beta_j,k_j,D_j>0,~~j=1,\ldots,N,
\end{align}
where $W_j$'s denote mutually independent (generalized) derivatives of the Wiener process \cite[e.g.]{sobczyk2013stochastic}. Each oscillator admits a stationary density given by
\begin{align}\label{eq_stdensity}
\rho_j(q=w_1,\dot{q}=w_2)= c_j\exp \left\{ \frac{-\beta_j}{D_j}\big(\frac{w_2^2}{2}+k_j\frac{w_1^2}{2}\big)\right\} ,~~j=1,\ldots,N.
\end{align}  
with $c_j$ being the normalization constant \cite{sobczyk2013stochastic}. By setting $D_j=k_j\beta_j$, we make sure that the displacement of each oscillator has a standard normal distribution.  Next, in order to make the model replicate the \emph{dynamics} of the time series,  we optimize the free parameters of each system to match the PSD of $q$. To be more precise, let $S_j$ denote the power spectral density (PSD) of the stochastic process $\{U^tq_j\}$. On the other hand, the PSD of response, $q_j$, in \eqref{eq_SDE} is given as
\begin{align}\label{eq_WKrelation} 
\tilde{S}_j(\omega)= \frac{2D_j}{\big| k_j+i\beta_j\omega-\omega^2\big|^2},
\end{align}  
where $\omega$ denotes the frequency \cite{sobczyk2013stochastic}. We find the pair $(k_j,\beta_j)$ that minimizes the spectral difference
\begin{align}\label{eq_specdiff} 
\Delta= \int_{0}^{\frac{\omega_s}{2}} \big| \tilde{S}_j(\omega)-S_j(\omega) \big| d\omega
\end{align}  
with $\omega_s$ being the sampling frequency. We approximate $S_j(\omega)$ from the time series of $q_j=T_j(y)$ using the Welch method \cite{welch1967use}, and solve this 2D optimization problem using particle swarm optimization \footnote{Other methods, including gradient-based descent, can be used for this problem as well, but we have empirically found the particle swarm method to be very efficient in this setup}. Note that the spectrum of $y_j$  will be affected by the spectrum of $q_1,\ldots,q_j$ due to the coupling of variables in $T^{-1}$, however, we have found that minimizing the error of $q$'s spectra independently is an efficient way to minimize the error in the spectrum of $y$ --- see  \cref{fig_cavityspectra} for a 10-dimensional example. When the time series are not long enough for a robust approximation of PSD, we recommend matching the autocorrelation function up to some finite time lag $\tau$. As discussed in \cite{korda2018data}, the autocorrelation sequence corresponds to the moments of the  spectral measure and matching those moments in the limit of $\tau\rar\infty$ leads to matching of the PSD. 
At the end, the data-driven stochastic model for the system in \eqref{eq_dynamics} consists of the stochastic oscillators in \eqref{eq_SDE} with the observation map $T^{-1}$.

The crucial finding of our study is that this choice of hypothesis space, i.e. stochastic linear oscillators observed under the inverse of triangular polynomial maps of order two or three,  
provides an \changed{adjustable} trade-off between learning capacity and data-efficiency in learning models of modal dynamics in strongly nonlinear flows with non-Gaussian statistics.  This approach is also advantageous since it is easily scalable to tens of dimensions: for example, the transport maps of the cavity flow example in \cref{sec_cavity}, involving 25000 data points in 10 dimensions, are computed in a few minutes on a plain desktop. The potential drawback is that our hypothesis space is relatively small (e.g., compared to deep neural networks) and there is no guarantee that  it would be suitable for learning other classes of nonlinear systems that may arise in other disciplines. In principle, however,  one could use a more sophisticated class of reference measures (e.g. bimodal distributions) and stochastic models  (e.g. non-Markovian linear stochastic systems), or even a more inclusive space of transport maps,  based on the complexity of the target system and availability of data.

\begin{figure}[!h]
        \centerline{\includegraphics[width=1 \textwidth]{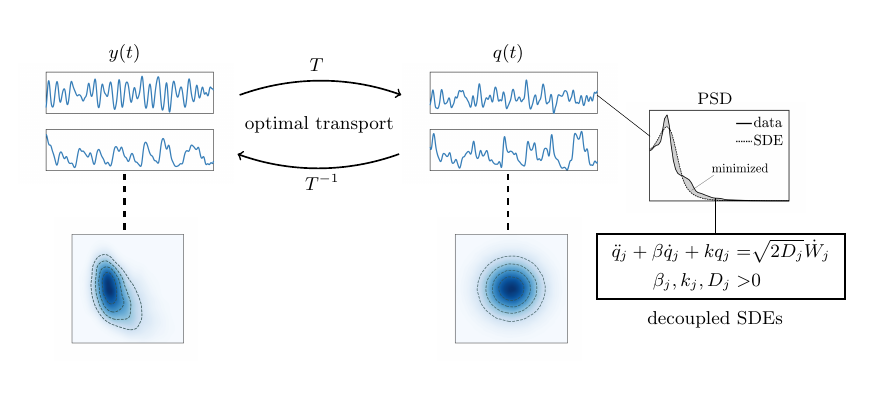}}
        \caption{\footnotesize \textbf{Framework for stochastic modeling.} A change of variable is found that transports the data distribution to standard normal distribution and then a stochastic oscillator is built by minimizing the difference between the PSD of transported time series and the SDE model (the shaded area).}
        \label{fig_framework} 
\end{figure}

Finally, we note that our framework for modeling, summarized in \cref{fig_framework},  is  connected to other methods of nonlinear systems identification: 
from the systems and signals perspective, our approach resembles the classical Wiener system identification which approximates nonlinear systems as a linear system with a nonlinear observation map  \cite{rugh1981nonlinear,boyd1985fading}. \changed{Our approach can also be considered as a data-driven discovery of a conjugacy \cite{wiggins2003introduction} between the data coordinates and the space of linear stochastic oscillators.  Our conjugacy map (the transport map) is particularly designed so that the obtained conjugate system has statistically independent coordinates.   Our approach is similar in spirit to other data-driven methods for discovering conjugate systems with independent components of motion \cite{mezic2017koopman} (i.e. phases along each direction of the tori or equivalently principal Koopman  eigenfunctions ). 
}
 Finally, our method provides a backward technique for fitting solutions of Fokker-Planck equations to data: direct construction of an SDE which gives rise to the data distribution requires solving Fokker-Planck equation in large dimensions which is often difficult or impossible, but in this framework the data is mapped to some distribution which is already a solution to a well-known set of models.   
 

\clearpage
\section{Results} 
\subsection{Lorenz 96 system}\label{sec_lorenz}
\begin{figure}[!h]
        \centerline{\includegraphics[width=1 \textwidth]{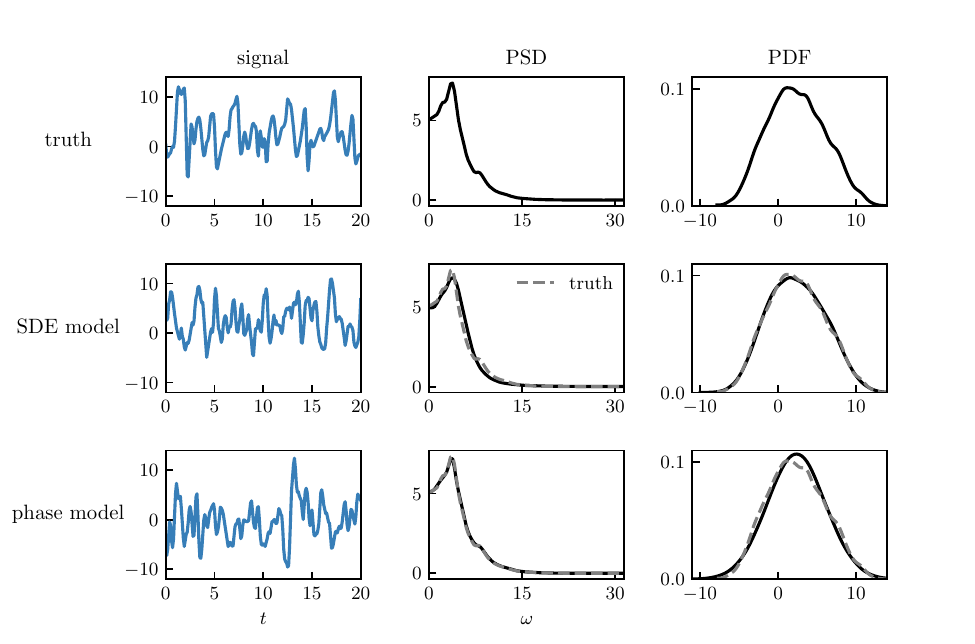}}
        \caption{\footnotesize \textbf{Data-driven modeling of Lorenz 96 system.} Top: truth observation on a state observable of chaotic Lorenz 96 system, middle: approximation via our stochastic modeling framework, bottom: random phase model approximation.  The SDE model captures the skewness in the PDF and closely matches the spectrum. The quasi-periodic phase model, by design, matches the spectrum of data but systematically misses the skewness in the PDF. The PSDs are computed after removing the mean. }
        \label{fig_lorenz}
\end{figure}

Lorenz 96 system  is a toy model of atmospheric turbulence and is commonly used in assessment of modeling, closure and data assimilation schemes for strongly nonlinear systems \cite{majda2005information}.  This model describes a set of interacting states on a ring which evolve according to 
\begin{align}\label{eq_lorenz}
\dot{x}_k =& x_{k-1}(x_{k+1}-x_{k-2})-x_k+F,\qquad k=1,2,\ldots,K\\
&x_0=x_K,\quad x_{K+1}=x_1,\quad x_{-1}=x_{K-1}. \notag
\end{align}
 For the standard parameter values of $K=40$ and $F=8$, trajectories converge to a chaotic attractor. We assume that we can only observe the first component of the state while the trajectory is evolving on the attractor, i.e.,
\begin{align}\label{eq_lorenz}
y = g(\x)= x_1.
\end{align} 
We record a time series of $y$ with the length of $1000$ and sampling interval of 0.1. After applying our SDE modeling framework to this time series, we find the stochastic model to be
\begin{align}\label{eq_lorenz_sde}
\ddot{q} + 4.73 \dot{q} +26.26 q= 15.76 \dot{W},
\end{align} 
where the state $q$ is related to the observable $y$ through the map \footnote{In the actual computation, the polynomial is not realized as in \eqref{eq_lorenz_tm}. See \cref{sec_TMmodeling} for details.}
\begin{align}\label{eq_lorenz_tm}
q=T(y)\approx 0.001y^3-0.010y^2+0.279y-0.570.  
\end{align}
The probability density function (PDF) and PSD of the observable $y=T^{-1}(q)$ are shown in  \cref{fig_lorenz}.
Due to the inclusion of nonlinear terms in the transport, this model is capable of capturing the skewness in the PDF of observable $y$. This is a crucial feature for models of turbulent systems since skewness gives rise to the nonlinear cascade of energy, and (linear) models with Gaussian distribution are incapable of describing this type of energy transfer between the system components.

We compare our model  to a data-driven approximation of the Koopman operator for the Lorenz system. The Koopman operator is usually approximated via  Extended Dynamic Mode Decomposition (EDMD) algorithm \cite{williams2015data,klus2015numerical,korda2017convergence}, however, the choice of observable dictionary in EDMD, which plays a critical role in the approximation, is not generally well-resolved. The most popular choice is to include delay-embeddings of the measurements \cite{arbabi2017ergodic}, but it is known that for measure-preserving chaotic systems the finite-dimensional approximations possess only eigenvalues with negative real part which lead to trivial statistics \cite{korda2018data}.
Here, we use the phase model approximation of the observable evolution in the form of 
\begin{align}\label{eq_rpm} 
y(t) = \sum_{j=1}^n \alpha_{j,y} e^{i(\omega_j t + \zeta_j)}
\end{align} 
where $\omega_j$'s are randomly chosen frequencies from the support of the Koopman spectral measure for $y$, $\alpha_j$'s are determined from spectral projections of $y$, and $\zeta_j$'s are initial phases randomly drawn from $\U[0,2\pi)$ (see Appendix \ref{sec_rpm} for more detail). It is shown in \cite{korda2018data} that, as $n\rar\infty$, such an approximation converges to the true evolution of  $y$ over \emph{finite time intervals}.  A realization of this quasi-periodic phase model for the Lorenz state for $n=200$ is shown in \cref{fig_lorenz}. This model, by construction, accurately reproduces the spectrum of the data  but fails to capture the skewness in the PDF. In fact, the phase model, which is commonly used for ocean wave modeling, is shown to converge to Gaussian statistics as $n\rar \infty$ \cite{kinsman1984wind}. This example highlights the fact that the current data-driven approximations of the Koopman operator in the form of finite-dimensional linear systems cannot reproduce complicated statistics observed in measure-preserving chaotic systems. \changed{We note that there have been a few recent works that aim to find nonlinear or stochastic representation of Koopman operator that do not have limitations of the quasi-periodic models \cite{brunton2016chaos,mezic2020spectrum} but the application of those approaches to strongly nonlinear flows remains to be explored.}

\clearpage
\subsection{Chaotic cavity flow and Spectral Proper Orthogonal Decomposition (SPOD)} \label{sec_cavity}
Fluid flows at large length and velocity scales are canonical examples of high-dimensional chaotic behavior. We discuss the application of our framework to such flows using a numerical model of the chaotic lid-driven cavity flow. This flow consists of an incompressible fluid in a 2D square domain, $\mathcal{D}=[-1,1]^2$, with solid walls, where the steady sliding of the top wall, given by
 \begin{equation}\label{eq_toplid}
u(y=1)=(1-x^2)^2,\quad -1\leq x \leq 1,
\end{equation} 
induces a circulatory fluid motion inside the cavity. The Reynolds number of this flow is defined as $Re=2/\kappa$ where $\kappa$ is the kinematic viscosity of the fluid in the numerical simulation. 
The cavity flow at $Re=30,000$ converges to a measure-preserving chaotic attractor exhibiting a purely continuous Koopman  spectrum \cite{arbabi2017study}. Here, we model the evolution of the post-transient cavity flow in two steps: first, we use the modes obtained by Spectral Proper Orthogonal Decomposition (SPOD) \cite{lumley1970stochastic,towne2018spectral} as a spatial basis for description of the flow evolution. We justify the use of SPOD modes through its connection with the spectral expansion  of the Koopman operator with continuous spectrum. In the second step, we use the framework based on measure transport and spectral matching to find the stochastic model that captures the evolution of flow in the SPOD coordinates.

Consider a non-homogenous stationary turbulent flow.  The velocity field at each point is a random variable defined on the underlying measure-preserving attractor, but due to non-homogeneity, velocity at different points have different statistics and spectra. However, since all these variables arise from the same underlying attractor, we can define characteristic spatial fields that connect the spectra of various variables through the notion of Koopman spectral measure. Let  $\ub_\x$ and $\ub_{\x'}$ denote velocity field at location $\x$ and $\x'$ in the flow domain. The spectral expansion of the Koopman operator for these two random variables is
 \begin{equation}\label{eq_KoopmanCross}
<\ub(\x),U^\tau \ub(\x')>_\mu = \int_{0}^{2\pi}e^{ i \omega \tau} \rho_{\x,\x'}(\omega)d\omega,
\end{equation} 
where $ \rho_{\x,\x'}$ is the Koopman cross spectral density of $\ub_\x$ and $\ub_{\x'}$, and $<\cdot,\cdot>_\mu$ is the inner product with respect to the invariant measure on the attractor \cite{mezic2005}. The SPOD modes of the flow at frequency $\omega$, denoted by $\psi(\x,\omega)$,  are defined as solutions of the eigenvalue problem \footnote{In contrast to \cite{towne2018spectral} we do not include the time-dependent oscillation in the definition of SPOD modes, since for chaotic systems there are no observables with single-frequency oscillations. }
\begin{equation}
  \int_{\Omega}\rho_{\x,\x'}(\omega)\psi(\mathbf{x}',\omega)d\mathbf{x}'
  =\lambda\psi(\mathbf{x},\omega).
\end{equation}
As such, SPOD modes are intrinsic dynamical properties that do not depend on the choice of flow realization, as opposed to dynamic modes or Fourier modes, and therefore provide a robust choice of basis for description of the flow evolution. A detailed discussion of Koopman spectral density and its connection with SPOD is given in Appendix \ref{sec_koopmanspec}.

We compute the SPOD modes of cavity flow using the algorithm in \cite{towne2018spectral}. We define the SPOD coordinates to be the projection of the flow onto the 10 most energetic SPOD modes at different frequencies, which contain $\sim 50\%$ of the turbulent kinetic energy, i.e.,
 \begin{align}  \label{eq_SPODy}
 y_j(t)=<\mathbf{u}(\x,t),\psi_j>_{\mathcal{D}}, \qquad j=1,2,\ldots,10.
 \end{align}
 where $<\cdot,\cdot>_\mathcal{D}$ denotes the spatial inner product over the flow domain.

\begin{figure}[H]
\centering
\subfloat[][]{\includegraphics[width=3.2in]{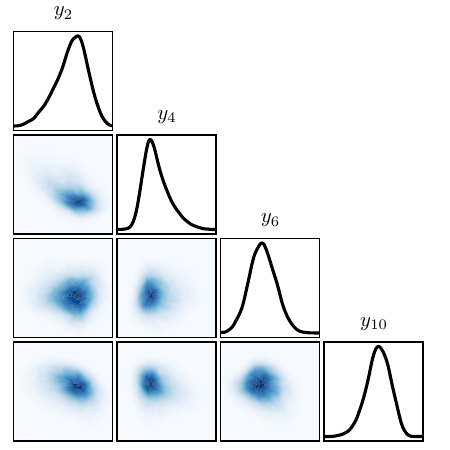}} ~
\subfloat[][]{\includegraphics[width=3.2in]{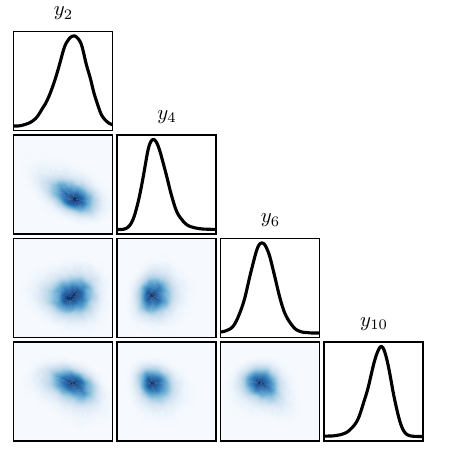}}

\subfloat[][]{\includegraphics[width=6.5in]{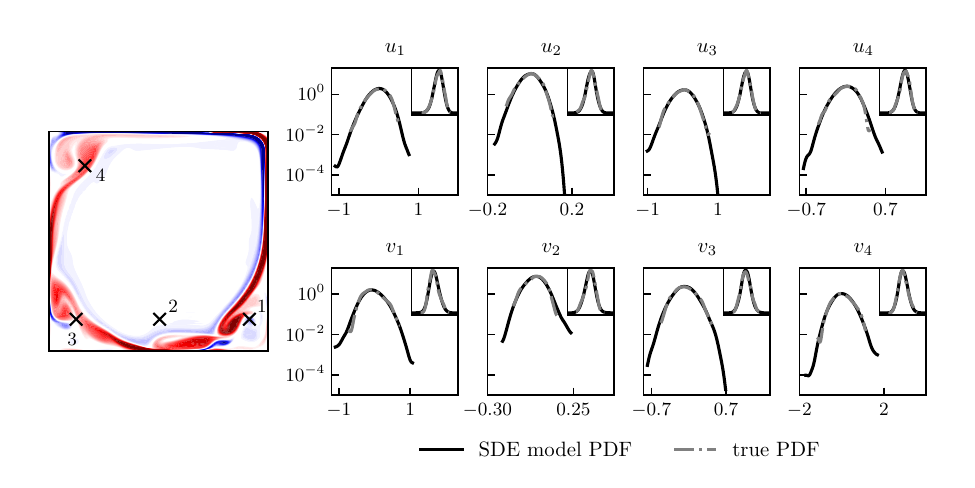}}
\caption{\textbf{Modeling of chaotic cavity flow at $\mathbf{Re=30000}$.} a) Single and pairwise marginal PDFs of 4 (out of 10) SPOD coordinates from data, b) same marginals by the SDE model. The quantile axes limits are  $(-0.026,0.026)$. c) Location of sensors for velocity measurements and a vorticity snapshot of the cavity flow (left) and the log PDFs of $u-$ and $v-$velocity generated by the SDE model and the 10-dimensional representation of the flow in SPOD coordinates (labeled as truth). The insets show the PDFs on a linear ordinate scale. The matching between the pointwise statistics indicates the skill of the SDE model in capturing the 10-dimensional joint pdf of SPOD coordinates. }
\label{fig_cavity_modeling}
\end{figure}

As the training data for our SDE model, we use a single time-series of SPOD coordinates with the length of 2500 seconds and sampling rate of 10 \emph{Hz}. We use a polynomial map of (total) degree 3 to compute the transport between the distribution of SPOD coordinates  and the standard normal distribution, and identify the corresponding SDE. After generating a 10000-second-long trajectory of the SDE, we compute the statistics of that trajectory, observed under the inverse of polynomial map. \Cref{fig_cavity_modeling} (a,b) shows the excellent agreement between the true marginals and the ones obtained from the SDE model. 
Next, we use the SPOD data generated by the model to construct a new flow trajectory and compare the statistics of pointwise velocity measurements with the original data projected to the 10-dimensional subspace of SPOD modes. The results, shown in \cref{fig_cavity_modeling}(c), confirm the accurate recovery of pointwise statistics hence indicating the high skill of our model in capturing the 10-dimensional joint pdf of flow modal data. Moreover, the match between log PDFs shows the utility of polynomial transport for data-driven extrapolation of PDF tails which we will explore in the next example.  
More details on modeling of cavity flow including the structure of SPOD modes and all 2d marginals can be found in Appendix \ref{sec_cavityextra}.

\clearpage
\subsection{Reanalysis climate data and tail extrapolation}
Statistical characterization of extreme events, i.e. probabilistically rare events that arise from combination of dynamics and randomness, is an important and challenging task. The importance stems from the significant and disruptive effect of extreme events such as earthquakes, rogue waves  and extreme weather patterns. Yet finding probabilities of extreme events, which reside at the tails of the PDF, is challenging because it requires a very large number of data points. The theoretical setup for study of such events is the extreme value theory (e.g. \cite{coles2001introduction}) which has been recently extended to chaotic dynamical system with some success \cite{freitas2008link,lucarini2016extremes}, but remains limited to systems with specific mixing properties. In the recent years, there has been a number of methods that use model reduction or machine learning  for data-efficient quantification of extreme events \cite{mohamad2016probabilistic,mohamad2018sequential}. 

Our modeling framework can be used for characterization of extreme events using relatively little data.
The key idea is that we can model a large space of probability measures with heavy tails by using transport to a reference measure like standard normal distribution
By discovering such a transport map, we identify our data distribution with a pullback of Gaussian probability under that map and the resulting distribution provides an approximation of the tails of the real data distribution. Although analytic description of the tails in these distributions are not available when the transport map is a high-degree multivariate polynomial map, we can easily generate a large sample from the reference measure and pass them through the transport map inverse to construct the tail approximation. In the following, we consider the application of this method to obtain a short extrapolation of super-Gaussian tails in climate data. 

\begin{figure}[!h]
        \centerline{\includegraphics[width=1 \textwidth]{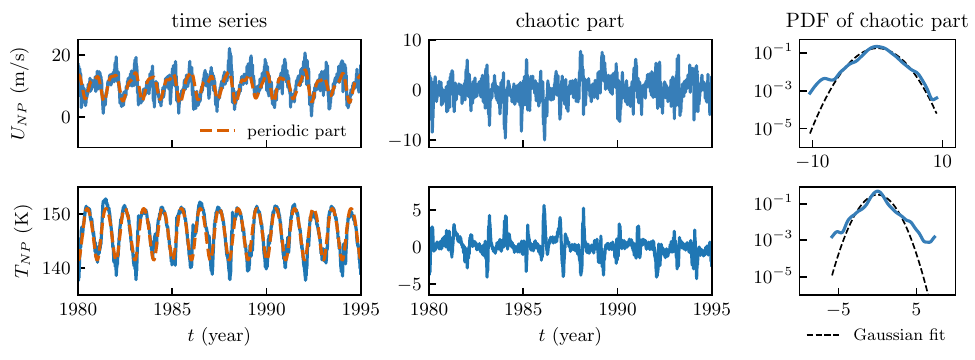}}
        \caption{\footnotesize \textbf{Velocity and temperature at the North Pole:} The time series showing wavelet coefficients of $u$-velocity  and temperature field at North Pole (left). We extract the periodic component consisting of Fourier modes with at least 1\% of signal variance to isolate the chaotic part (middle) which shows super-Gaussian tails (right). }
        \label{fig_cl1} 
\end{figure}

Our data is based on 6-hourly reanalysis of velocity and temperature in the earth atmosphere recorded at the 100-\emph{mbar} iso-pressure surface from 1981 to 2017 \cite{berrisford}.  These global fields are then expanded in a spherical wavelet basis described in \cite{lessig}. \Cref{fig_cl1} shows the time series  of  the expansion coefficients for a level-1 wavelet envelope centered on top of the North Pole. These time series exhibit a combination of periodic and chaotic oscillations. Therefore, we model them as observation on a stationary system with mixed spectra, i.e., a system that possesses both discrete and continuous Koopman spectra (see e.g. \cite{broer1993mixed,mezic2005}). In the first step, we extract the periodic component of data (i.e. Koopman modes) by removing Fourier modes of time series that correspond to more than 1\% of fluctuation energy. We apply our stochastic modeling framework to the chaotic remainders with heavy tails.

\begin{figure}[H]
\centering
\subfloat[][]{\includegraphics[width=6.5in]{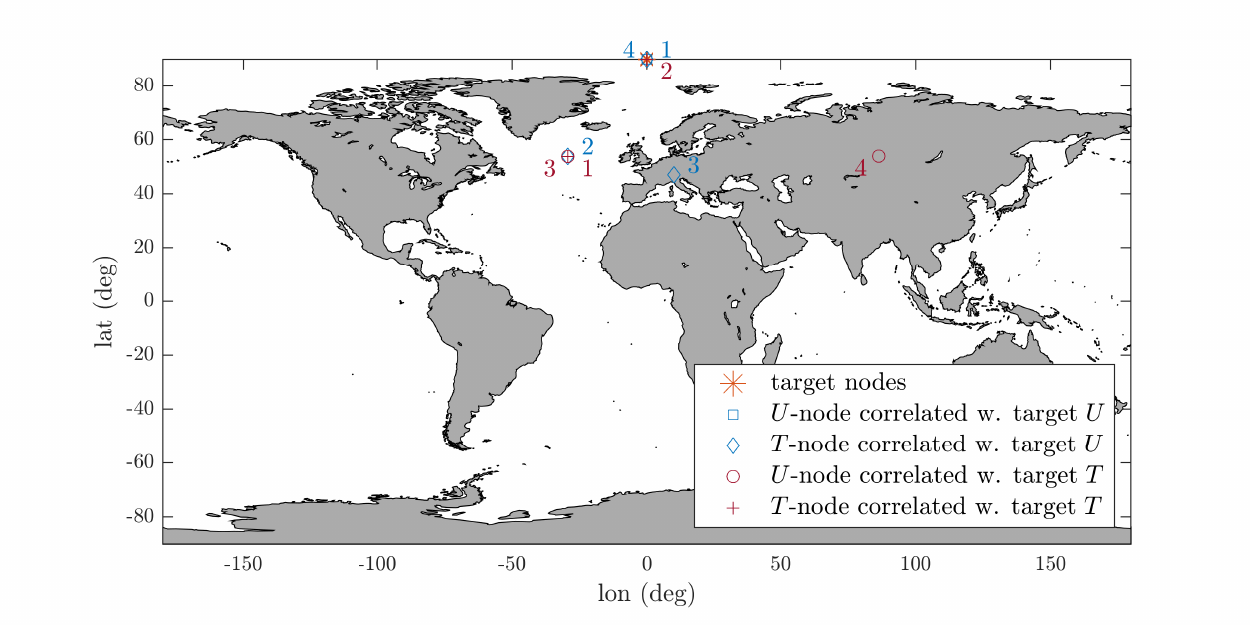}}

\subfloat[][]{\includegraphics[width=6.5in]{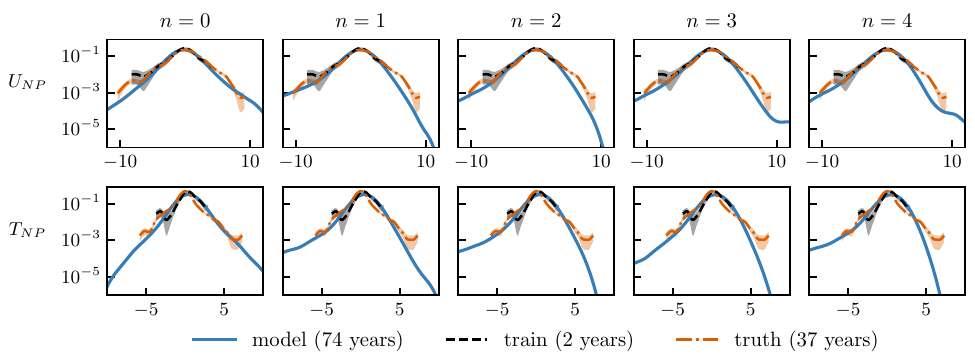}}
\caption{\textbf{Extrapolation of tails for climate data.} a) Location of covariates used to build models for the target random variables, i.e., $u$-velocity and temperature at the North Pole, b) approximation of PDF tails by generating a surrogate trajectory of the SDE model and pulling it back under the transport map. $n$ is the number of covariates  used to learn the SDE model. The training data is the time series in 1981-1982 and the truth data is the time series in 1981-2017. The shaded envelopes show the 95\% (pointwise) confidence interval of PDF estimation for training and truth data.
}
\label{fig_cl2}
\end{figure}

We construct a sequence of stochastic models for each of the chaotic time series obtained from measurements at the North Pole. In each model, we add an extra random variable from velocity and temperature measurements at other locations to capture the statistical dependence between those variables and the target variable (i.e. chaotic part of temperature or $u$-velocity). These covariates are chosen from the set of all wavelet coefficients ($>8000$ variables) in the order of decreasing absolute linear correlation with the target variables. The geographic position of the covariates is shown in \cref{fig_cl2}(a). 
Importantly, we exploit the triangular structure of transport in \eqref{eq_triangular} and place the target variable on the bottom of each map, so that its mapping to the reference measure is informed by all the covariates. 
For training the models we use only the time series from 1981 and 1982 and polynomials of total degree 2.

After discovering each model, we generate a 74-year long trajectory of the SDE model and pull it back to the space of observations. 
The PDFs of the observed map provide a long-time extrapolation to the PDFs of the training data. As shown in \cref{fig_cl2}, those PDFs are in good agreement with the original reanalysis data of 37 years. In particular, the approximation of tails with good training data (i.e. left tails in the top and bottom row) show consistency with the addition of covariates  to the model.

\clearpage
\section{Dependence on training data and degree of polynomial mapping}
Under our assumptions for modeling, the solution to the transport map is known to exist, however, restricting the search to a specific set of functions, such as polynomials, induces a bias in approximating the statistics of the flow.  We are usually interested in evaluating expected value of a square-integrable random variable $h$, and it is known \cite{marzouk2016sampling} that 
 \begin{align}  \label{eq_error_KL}
 \|\Ex_\nu[h]-\Ex_{\tilde\nu}[h]\| = \sqrt{2 (\Ex_\nu[\|h\|^2]-\Ex_{\tilde\nu}[\|h\|^2]} \sqrt{ \Ex_{\nu}\bigg[\log \frac{\nu}{\tilde \nu} \bigg]}
  \end{align}
where $\nu$ is the target distribution and $\tilde\nu$ is the distribution generated by our model. Note that the last term on the right-hand-side is the objective of the optimization problem in \cref{sec_TMmodeling} and therefore we can control the bias in estimation of the statistics through our approximation. Rigorous results on convergence for the type of transport maps we use here are sparse (for example, see \cite{zech2020sparse} for a similar setup bounded domains), and theoretical characterization of  the bias in application to strongly nonlinear flows is beyond the scope of this paper. Instead, we perform a numerical study of convergence for the examples of Lorenz system and cavity flow where we have access to extensive data as the ground truth.

\begin{figure}[H]
\centering
\subfloat[][]{\includegraphics[scale=1]{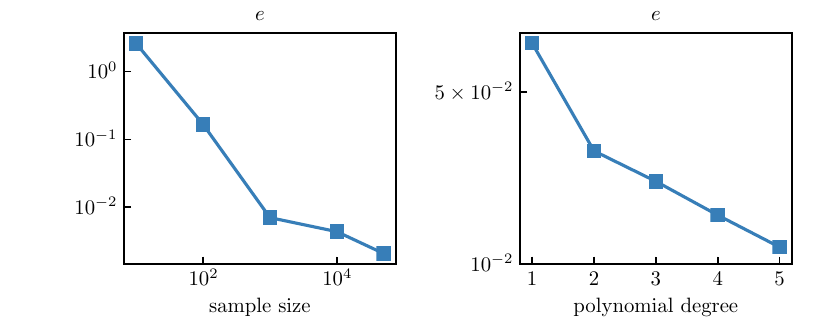}}

\subfloat[][]{\includegraphics[scale=1]{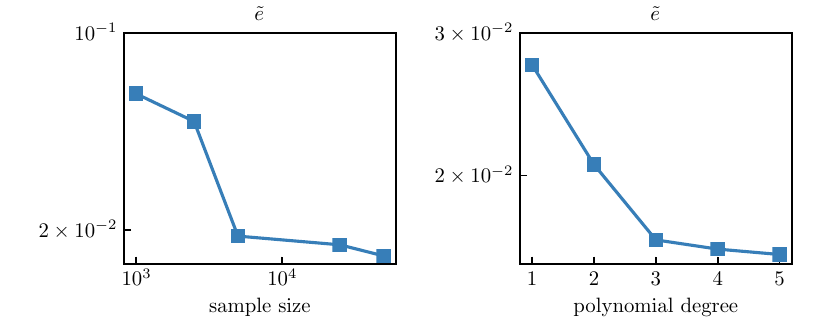}}
\caption{\textbf{Statistical error for modeling via transport} a) Lorenz 96 system: the PDF variance error of state variable (left)  versus sample size with fixed polynomial degree 3  and (right) versus polynomial degree with fixed sample size of 1000. b) Cavity flow: the average of PDF variance error for 
pointwise velocity  (right) versus sample size with fixed polynomial degree 2 and  (right) versus polynomial degree with fixed sample size 25000. }
\label{fig_analysis} 
\end{figure}

To assess the quality of modeling via transport we use  the  \emph{the variance diagnostic} \cite{marzouk2016sampling},
 \begin{align}  \label{eq_var_error}
 e=\mathbb{V}\text{ar}_{\rho_\text{truth}} \bigg( \log \frac{\rho_\text{truth}}{\rho_\text{model}} \bigg).
 \end{align}
 where $\rho_\text{truth}$ and $\rho_\text{model}$  denote the density of the data distribution and the density of the distribution generated by the transport model.
In the case of Lorenz 96, we have direct access to the PDF of the state variable (computed from simulations and taken as truth) and we can evaluate the diagnostic directly. As shown in \cref{fig_analysis}(a), with the increase of the training sample size and the polynomial degree of the transport map, the variance diagnostic rapidly decreases, indicating the scalable accuracy of the transport-based model for Lorenz 96 system.

In the cavity flow example, the target distribution is given by sample data in a 10-dimensional  space and there is no analytic expression for its density. Therefore, as a proxy for the error of modeling, we look at the variance diagnostic in recovering the statistics of a set of observables of the flow. Specifically, we use the statistics of pointwise velocity measurements because they are functions of all the SPOD coordinates, and   therefore reveal the performance of the model in emulating the \textit{high-dimensional joint} distribution of the data.
Let $e_{u(\x)},~e_{v(\x)}$ be the variance diagnostic for the pointwise PDF of, respectively, $u-$ and $v-$ velocity field at position $\x$ in the flow domain.  We compute the average of this error over 50 randomly placed sensors in the flow to assess the overall quality of this model:
 \begin{align}  \label{eq_var_error}
 \tilde e=\frac{1}{n}\sum_{k=1}^{50}e_{u(\x_k)}+e_{v(\x_k)}
  \end{align}
 The results in \cref{fig_analysis}(b) shows similar results: with the increase of training sample size and the polynomial degree, the proxy variance diagnostics decreases. This shows that the space of models considered here are good candidates for learning models of strongly nonlinear flows from relatively little data.

%

\clearpage
\section{Discussion}
We presented a framework for generative modeling of strongly nonlinear flows in the form of decoupled stochastic oscillators with nonlinear observation maps. The key feature of our framework is the use of measure transport to model non-Gaussian invariant measures that arise from statistics of complex chaotic systems.
Application to the high-Reynolds cavity flow showed that models generated by our framework can accurately reproduce the 10-dimensional joint pdfs of modal coordinates and hence recover the pointwise statistics in the flow.  We also showed the promise of our framework in data-driven characterization of extreme events through an example of reanalysis climate data.

In the context of operator-theoretic approximation, our framework is closest in spirit to the approach in \cite{dellnitz1999approximation,dellnitz2002dynamical} where the Perron-Frobenius operator of the underlying dynamical system is approximated using a Markov chain for  transition between the computational cells in the state space (also see \cite{govindarajan2018approximation} for a similar Koopman operator approximation). Although this approach reproduces both the statistics and dynamical behavior, it is not computationally scalable to moderate and high-dimensional systems, and for systems like complex flows one needs to use extreme coarse-graining of the dynamics e.g. by  clustering data into a few discrete representative states  \cite{kaiser2014cluster,perez2013identification}.    
More recent data-driven approaches for generative modeling are the Generative Adversarial Network (GAN)\cite{goodfellow2014generative} and Variational Auto-Encoders (VAEs) \cite{kingma2013auto}. GANs and VAEs have produced exemplary results in emulating high-dimensional and rich data distributions, however, like other deep networks, they offer little interpretability and lack of an inherent dynamical character. Training GANs, in particular, often requires non-trivial techniques and interestingly one of the effective improvements on GAN training, the Wasserstein GAN \cite{arjovsky2017wasserstein},  is based on ideas from measure transport. The framework presented in this paper provides a balanced trade-off for data-driven modeling of strongly nonlinear systems that arise in turbulent flows. Through the use of a single-layer triangular polynomial maps, it offers much interpretability in assessing the role of different variables on the observed dynamics, and due to separability of the underlying optimization problem it is easily extendable to tens of dimensions.

One promising direction for extension of our framework is to include physics-based constraints in computation of transport maps and the underlying SDE models. As shown in previous works, incorporating physical information about  the target system into the structure of learning can substantially increase the data-efficiency of modeling. In the context of statistical modeling for turbulent flows, statistical laws and tail slopes predicted from the theory provide suitable constraints. In this work, we showed the effectiveness of our approach for learning strongly nonlinear fluid flows, but our framework also has great potential for data-driven modeling in other applications such as nonlinear optics, chemical reaction networks, molecular dynamics and biological systems. For these systems, it might be necessary to make other judicious choices for the reference measure and the type of transport maps to achieve an efficient learning process.

\clearpage
\section*{Acknoweldgements}
We are grateful to Dr. Boyko Dodov and AIR\ Worldwide for providing the reanalysis climate data, as well as Prof. Christian Lessig  who performed the projection of the data to a spherical wavelet basis. We also thank Profs. Youssef Marzouk, Igor Mezi\'c and Yannis Kevrekidis for instructive discussions and pointing out related references. H.A. is grateful to Dr. Antoine Blanchard for notes on usage of NEK5000 and Dr. Daniele Bigoni on usage of the computational package for transport. This research was supported in part by the ARO-MURI
grant  W911NF-17-1-0306. T.S. acknowledges support from MIT Sea Grant through Doherty Associate Professorship and AIR\ Worldwide.

\section*{Data and source code }

Python Implementation of our framework and the data for producing the figures in this paper is available at \url{https://github.com/arbabiha/StochasticModelingwData} .

\clearpage

\begin{appendices}

 
\section{Data and compilation of pdfs}   

The data for Lorenz 96 system is generated using direct integration of \cref{eq_lorenz} using 4th-order Runge-Kutta with time step $\Delta t=0.01$. The velocity field for the lid-driven cavity flow is computed by direct solution of Navier-Stokes equations  in  NEK5000 \cite{nek5000-web-page}. The cavity domain is divided into 25 elements in each direction and the solution within each element is represented with a polynomial of 7th degree. The stationary solution is recorded after 3500 transient simulation time units. The processed reanalysis climate data is provided by AIR\ Worldwide.
The PDFs reported for Lorenz, cavity flow and climate (truth and model) are computed using a 100-bin histogram and, for visualizations, smoothed by convolving with Gaussian kernel with standard deviation of 2 bin widths. The PDFs for the climate training data are computed using 50 bins. The confidence intervals for PDFs of climate are computed using binomial statistics and the adjusted Wald formula in \cite{agresti1998approximate}.

 A Python implementation of our framework and the  data for producing the figures in this paper can be found  at \url{https://github.com/arbabiha/StochasticModelingwData}.

\section{Cross-spectral density of Koopman operator and Spectral Proper Orthogonal Decomposition (SPOD)}\label{sec_koopmanspec}
In this section, we recall the Koopman spectral expansion for chaotic systems and define the cross-spectral density of Koopman operator which is used in the definition of SPOD in \cref{sec_cavity}.
Consider a deterministic dynamical system with an attractor $\Omega$  which supports a physical measure $\mu$, and $\mu(\Omega)=1$.  Let $f,~g\in \h:=L^2(\Omega,\mu)$ be observables of this system, with the usual inner product 
\begin{align}
<f,g>_\mu=\int_Afg^*d\mu.
\end{align}
Now recall the Koopman operator $U^t$ defined as $Uf=f\circ F^t$ where $F^t$ is the reversible flow of the dynamical system. $U^t$ is a \emph{unitary} operator on $\h$, that is $(U^t)^*=U^{-t}$. This implies that the spectrum of the Koopman operator lies on the imaginary axis in the complex plane, and the spectral expansion of the Koopman operator  \cite{mezic2005,mezic2013analysis}  is given as
\begin{align}\label{eq:KSE}
U^tf =\sum_{k=0}^{\infty}v_k\phi_ke^{i\omega_kt} + \int_{0}^{2\pi}e^{i\omega t}dE_\omega f,
\end{align}
where the countable sum is the Koopman mode decomposition of $f$ associated with the quasi-periodic part of the evolution, $i\omega_j$ is a Koopman eigenvalue (i.e. an element of discrete spectrum) associated with eigenfunction $\phi_k$, and $E$ is the spectral measure of the Koopman operator associated with the continuous part of the spectrum. To be more precise, $E$ is a measure on $[0,2\pi)$ that takes values in the space of projections on $H$. That is, for every measurable set $B\subset[0,2\pi)$, $E_B$ is a projection operator, and $E_Bf$ is projection of $f$ onto the eigen-subspace associated with the part of spectrum residing in $B$. We are interested in continuous part of the Koopman spectrum which corresponds to chaotic behavior, and therefore \emph{we assume that there are no quasi-periodic parts, including the part with zero frequency (i.e. mean of the observable) present in the evolution.}

The operator-valued measure $E$ is difficult to characterize using data. Instead, we can define a positive real-valued measure associated only with $f$ and defined as follows \cite{maccluer2008elementary}, 
\begin{align}
\mu_f(B)=<E(B)f,f>_\mu.
\end{align}
Then, we can rewrite the spectral expansion of $f$ as
\begin{align}
<U^tf ,f>_\mu=\int_{0}^{2\pi}e^{i\omega t}<dE_\omega f,f>=\int_{0}^{2\pi}e^{i\omega t}d\mu_f(\omega)
\end{align}
where $\rho_f$ is the Koopman spectral density of observable $f$. Similarly, we can define a measure for spectral correlation of two distinct observables. That is
 \begin{equation}
<E(B)f, g>_\mu=<f, E(B)g>_\mu =\mu_{f,g}(B)
\end{equation}
defines a finite complex-valued measure $\mu_{f,g}$ on $[0,2\pi)$. Under the assumption of absolute continuity for this measure, we can write
 \begin{equation}\label{eq:Koopmancrossspectral}
\mu_{f,g}(B)=\int_B \rho_{f,g}d\alpha,
\end{equation}
with $\rho_{f,g}$ being \emph{Koopman spectral density of observables $f$ and $g$.} Consequently, we can write the following expansion for  the dynamic evolution of the two observables 
 \begin{equation}\label{eq:KoopmanCross}
<f,U^\tau g>_\mu = \int_{0}^{2\pi}e^{ i \omega \tau} \rho_{f,g}(\omega)d\omega.
\end{equation} 
In view of the duality between measure-preserving deterministic systems and stationary stochastic processes \cite{doob1953stochastic}, the expansion in \eqref{eq:KoopmanCross} is the same as the spectral expansion for stochastic processes used in definition of SPOD \cite{lumley1970stochastic,towne2018spectral}.

\section{Construction of random phase model for chaotic systems}\label{sec_rpm}
Let $g:\Omega\rar\mathbb{R}$ be an observable on the measure-preserving dynamical system. Given the spectral density of $g$,   \cref{alg:RPM} generates a random phase model for the evolution of $g$ in time. The main idea is to approximate the spectral measure of $g$ using sum of delta functions, that is, modeling the evolution as rotation on tori. Using ergodicity, we can use a single realization to compute the PDF reported in \cref{sec_lorenz}.

To see the connection with the approximation of the Koopman operator note that the work in \cite{korda2018data} has shown that approximations such as
 \begin{equation}\label{eq_milanqp}
U_m^1 g := \sum_{j=1}^{m}e^{2\pi i \omega_j} E(B_j) g,
\end{equation} 
where $B_j$'s are a partition of $[0,2\pi)$ and $\omega_j \in B_j$, will converge to the true evolution $U^1g$ in $L^2(\Omega,\mu)$ in the limit of infinitely refined partition. This approximation is in the function space and the choice of $B_j,\omega_j$'s are not unique. \Cref{alg:RPM} generates a realization of this approximation, with $B_j=[e_{j-1},e_j)$, that is spectrally consistent, i.e.,
 \begin{equation}\label{eq_milanqp}
|a_j|^2 = \int_{B_j} \rho_g(\omega) d\omega.
\end{equation}

\begin{algorithm}
\caption{Construction of random phase model from spectral density}\label{alg:RPM}
\begin{algorithmic}[1]
\Require spectral density $\rho(\omega)$, number of intervals on the frequency domain $m$ 
\State Draw $m$ random value of cell edges, $e_j$,  from uniform distribution on $[0,2\pi)$.
\State Set $e_0=0$ and $e_{m+1}=2\pi$.
\State Sort the random cell edges to form the sequence $e_j,~j=0,1,\ldots,m+1$ with $e_{j+1}>e_j$.
\For{$j=1,\ldots,m$}
\State Let $\omega_j = \frac{1}{2}(e_j+e_{j-1}),~j=1,\ldots,m$.
\State Let $\Delta\omega_j=e_j-e_{j-1},~j=1,\ldots,m$.
\State Let 
\begin{align}
a_j = \sqrt{\rho(\omega_j)\Delta\omega_j }\approx \bigg( \int_{e_{j-1}}^{e_j}\rho(\omega) d\omega\bigg)^{1/2}.
\end{align}
\EndFor
\State Draw $\zeta_j$ randomly and independently for each $j$ from the uniform distribution on $[0,2\pi)$ and let
\begin{align}\label{eq:rpm2}
g_t = \sum_{j=1}^{m} a_ke^{i(\omega_jt+\zeta_j)}
\end{align}
\end{algorithmic}
\end{algorithm}

\clearpage
\section{Supplemental information and figures for cavity flow}\label{sec_cavityextra}
\begin{figure}[ht!]
\centering
\subfloat[][]{\includegraphics[width=3.2in]{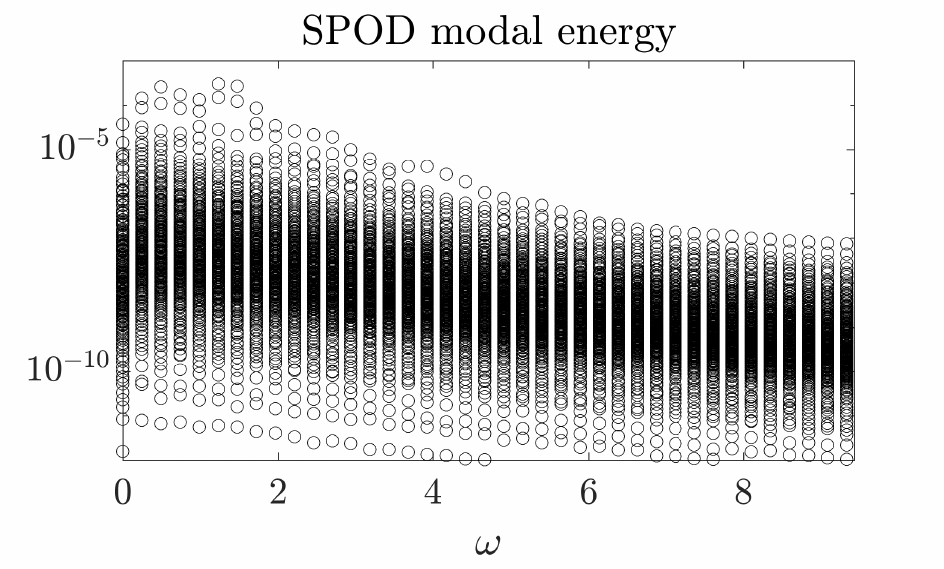}}

\subfloat[][]{\includegraphics[width=6.5in]{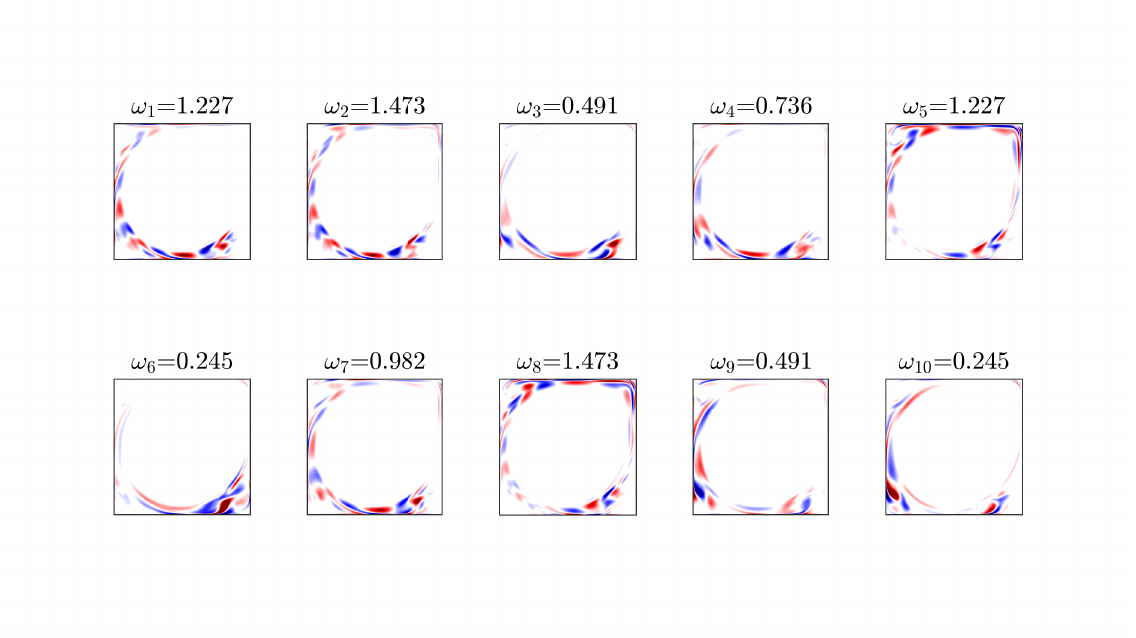}}
\caption{\textbf{Spectral orthogonal decomposition (SPOD) for cavity flow at $\mathbf{Re=30000}$.} a) Distribution of the kinetic energy within the SPOD modes vs frequency, for each frequency we have 155 spatially orthogonal modes. b) (real part of) vorticity for the top 10 energetic SPOD modes corresponding to the 10 highest circles in (a).  The red color marks clockwise rotation and the blue counter-clockwise.}
\label{fig_spod_spectrum}
\end{figure}

\begin{figure}[!h]
	\centerline{\includegraphics[width=1 \textwidth]{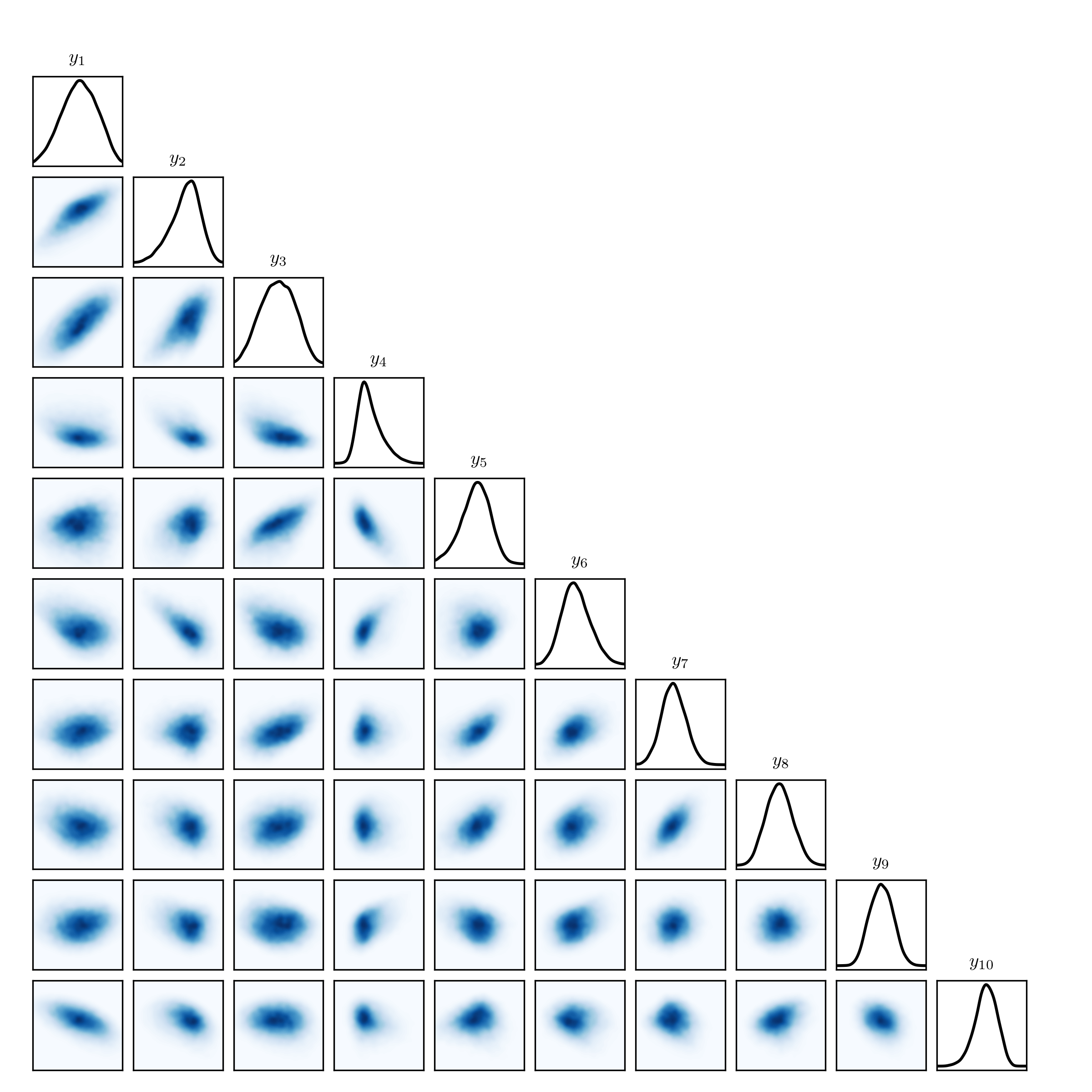}}
	\caption{\footnotesize  Marginal distributions of SPOD coordinates in cavity flow from the flow simulation. The quantile axes limits are  $(-0.026,0.026)$}
	\label{fig_cavity_10dmr}
\end{figure}

\begin{figure}[!h]
	\centerline{\includegraphics[width=1 \textwidth]{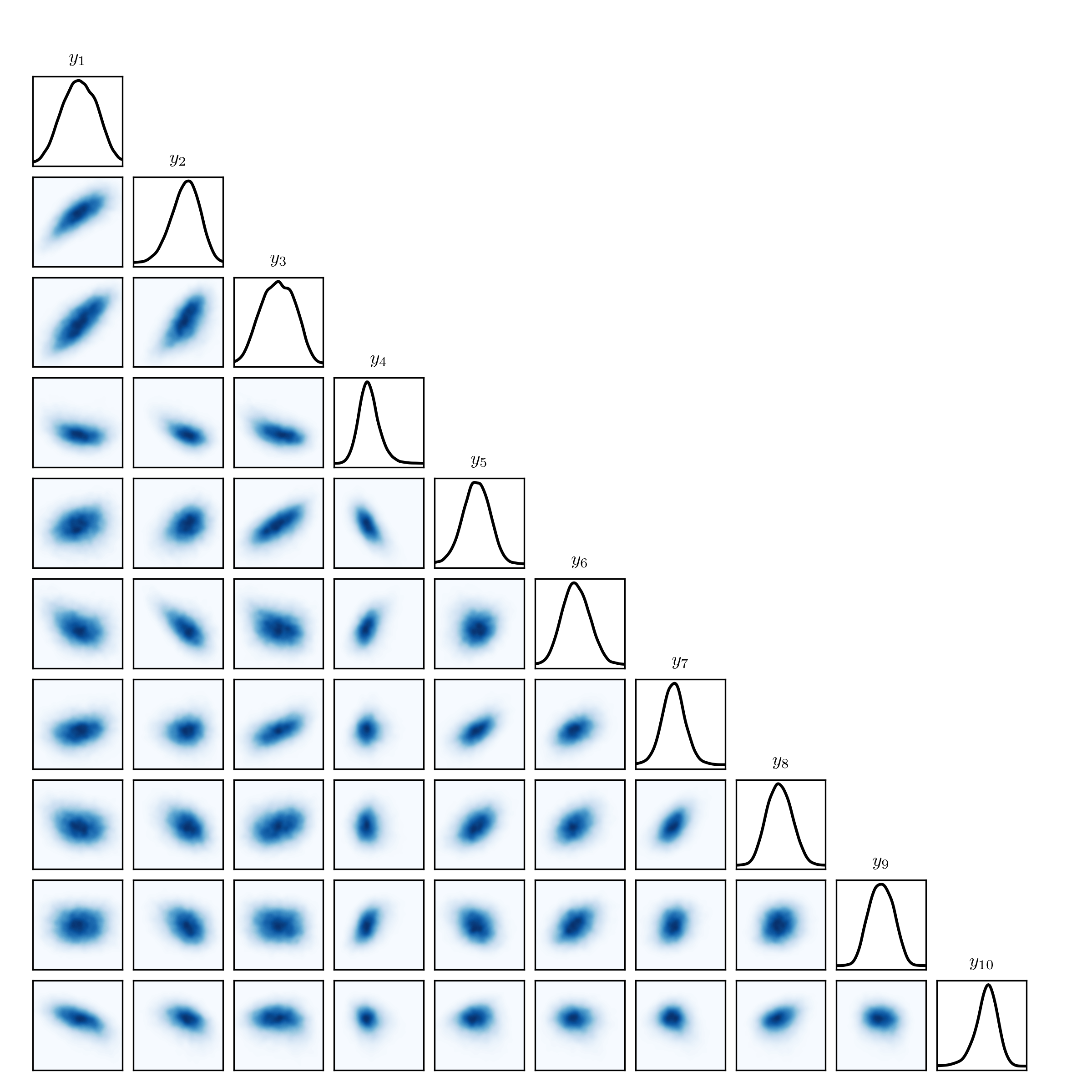}}
	\caption{\footnotesize  Marginal distributions of SPOD coordinates in cavity flow from the SDE model. The quantile axes limits are  $(-0.026,0.026)$. Compare to the truth in \cref{fig_cavity_10dmr}}
	\label{fig_cavity_10dmm}
\end{figure}

The data for the cavity flow consists of a single trajectory with the length of 12000 seconds. We have removed the mean flow from the data and then applied the algorithm in \cite{towne2018spectral} to the first 2000 seconds to compute the SPOD modes of flow shown in  \cref{fig_spod_spectrum}. Then projected the next 10000 seconds onto the top 10 energetic modes, using \cref{eq_SPODy}, to obtain the SPOD coordinates. We have used the first 2500 second of the SPOD coordinate time series, with sampling rate of 10 \emph{Hz}, for training of the SDE model.  The marginal distributions of flow data (\cref{fig_cavity_10dmr}) are computed using the whole trajectory. An SDE trajectory of the same length is used to compute the  marginal distributions in \cref{fig_cavity_modeling} and \cref{fig_cavity_10dmm}. To compute the pointwise statistics, we reconstruct the flow field via
\begin{align}
\tilde{\ub}(\x)=\sum_{j=1}^{N}\sum_{k=1}^{N} H_{jk} y_k \bs{\psi}_j(\x)
\end{align}
where the matrix $H$ is the inverse of Gramian matrix $G$ defined as $G_{jk}=<\bs{\psi}_j,\bs{\psi}_k>_\mathcal{D}$\footnote{The matrix $H$ appears to compensate the non-orthogonality of SPOD modes.}. The results in this paper are compiled using the real part of SPOD modes for simplicity in presentation.

\begin{figure}[!h]
	\centerline{\includegraphics[width=1 \textwidth]{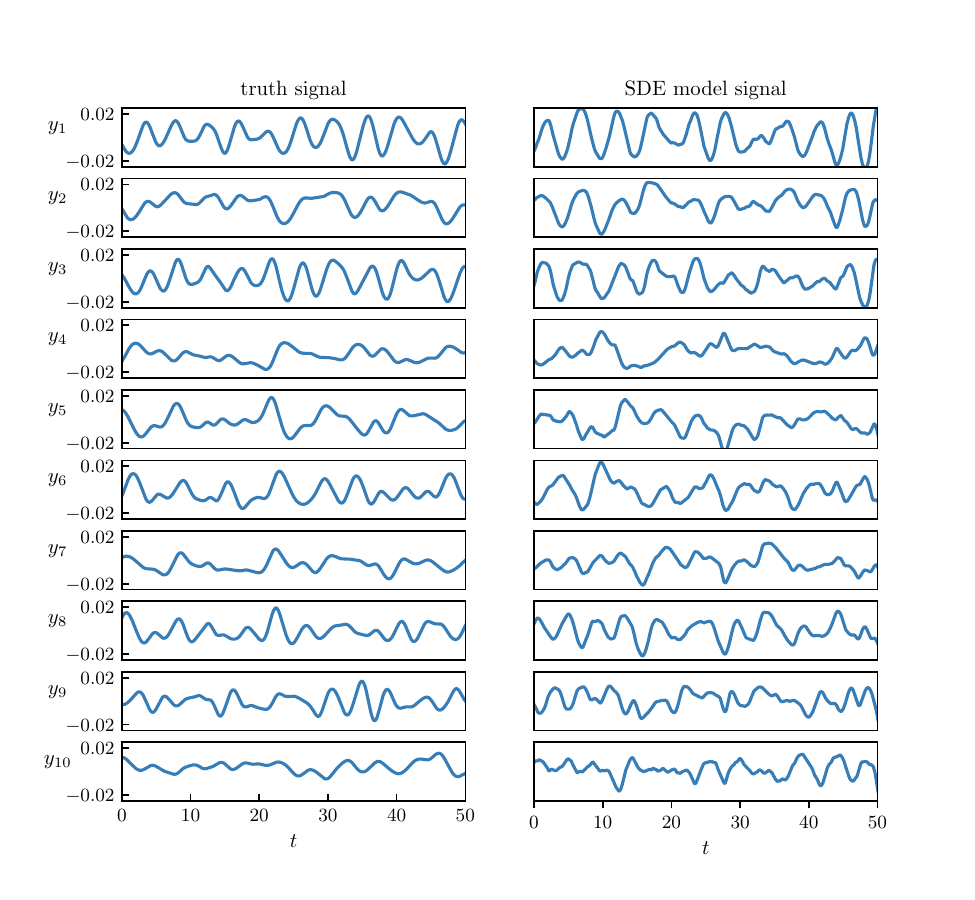}}
	\caption{\footnotesize  Time evolution of SPOD coordinates for cavity flow with samples from the truth data (numerical simulation of the flow, shown on the left) and a realization from the SDE model (right). }
	\label{fig_cavitysignal}
\end{figure}
 
 \begin{figure}[!h]
	\centerline{\includegraphics[width=1 \textwidth]{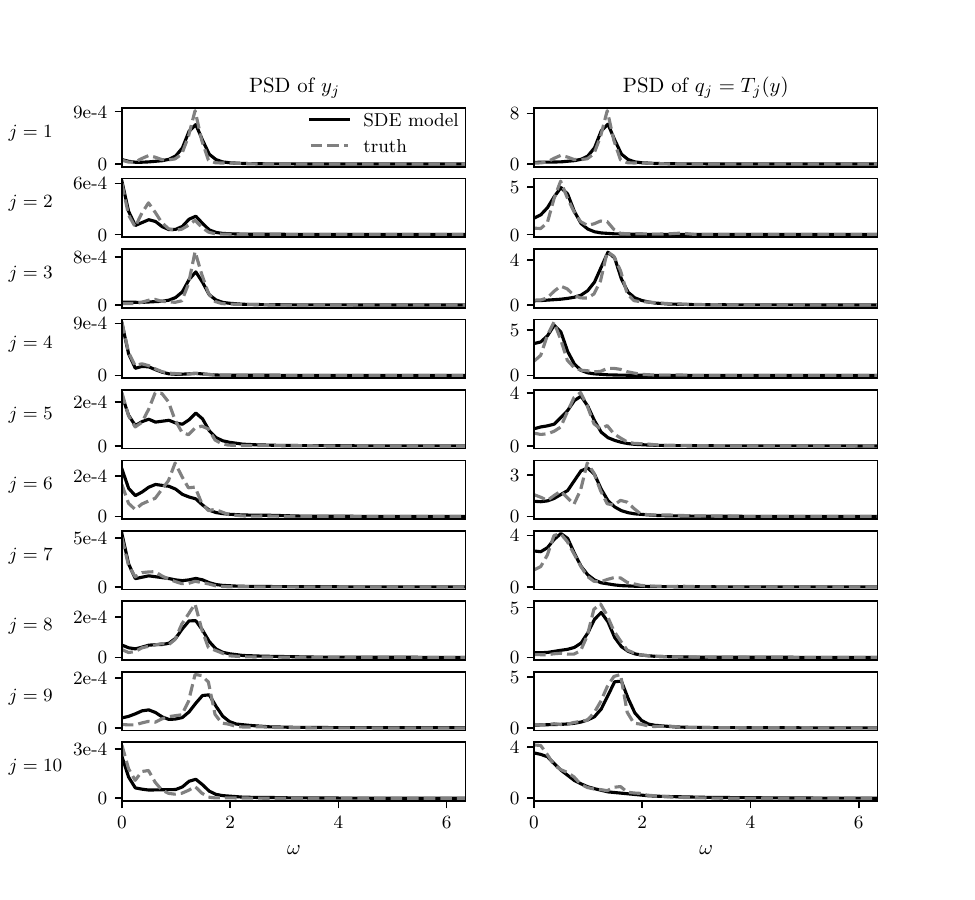}}
	\caption{\footnotesize  Power spectral density (PSD) of SPOD coordinates from the flow simulation and the SDE model (left), and the transported variable (right). }
	\label{fig_cavityspectra}
\end{figure}
 
 \end{appendices}

\clearpage
\bibliography{Stochastic_model_reduction.bbl}
\bibliographystyle{siam}
\end{document}